\documentclass[a4paper,12pt]{amsart}
\usepackage{amssymb}
\usepackage{ifthen}
\usepackage{graphicx}
\usepackage{color,xcolor}
\usepackage{array}
\usepackage{caption}
\nonstopmode
\numberwithin{equation}{section}
\newtheorem{thm}{Theorem}
\newtheorem{cor}{Corollary}
\newtheorem{lem}{Lemma}
\newtheorem{prop}{Proposition}

\newtheorem{conj}{Conjecture}
\newtheorem{prob}{Problem}

\theoremstyle{definition}
\newtheorem{defn}{Definition}
\newtheorem{ca}{Case}

\newtheorem{rem}{Remark}

\newenvironment{pf}[1][]{%
 \vskip 1mm
 \noindent
 \ifthenelse{\equal{#1}{}}%
  {{\slshape Proof. }}%
  {{\slshape #1.} }%
 }%
{\qed\medskip}

\newcounter{alphabet}


\newcounter{alphabet2}

\def\be{\begin{equation}}
\def\ee{\end{equation}}

\newcommand{\ben}{\begin{enumerate}}
\newcommand{\een}{\end{enumerate}}

\newcommand{\blem}{\begin{lem}}
\newcommand{\elem}{\end{lem}}
\newcommand{\bthm}{\begin{thm}}
\newcommand{\ethm}{\end{thm}}
\newcommand{\bcor}{\begin{cor}}
\newcommand{\ecor}{\end{cor}}
\newcommand{\beg}{\begin{exam}}
\newcommand{\eeg}{\end{exam}}
\newcommand{\begs}{\begin{examples}}
\newcommand{\eegs}{\end{examples}}
\newcommand{\bdefe}{\begin{defn}}
\newcommand{\edefe}{\end{defn}}
\newcommand{\bprob}{\begin{prob}}
\newcommand{\eprob}{\end{prob}}
\newcommand{\bques}{\begin{ques}}
\newcommand{\eques}{\end{ques}}
\newcommand{\bei}{\begin{itemize}}
\newcommand{\eei}{\end{itemize}}
\newcommand{\bcon}{\begin{conj}}
\newcommand{\econ}{\end{conj}}
\newcommand{\bop}{\begin{op}}
\newcommand{\eop}{\end{op}}

\newcommand{\bas}{\begin{assertion}}
\newcommand{\eas}{\end{assertion}}

\newcommand{\bfa}{\begin{fact}}
\newcommand{\efa}{\end{fact}}

\newcommand{\bca}{\begin{ca}}
\newcommand{\eca}{\end{ca}}

\newcommand{\bst}{\begin{step}}
\newcommand{\est}{\end{step}}

\newcommand{\bsca}{\begin{sca}}
\newcommand{\esca}{\end{sca}}

\newcommand{\bcl}{\begin{cl}}
\newcommand{\ecl}{\end{cl}}

\newcommand{\bmlem}{\begin{mlem}}
\newcommand{\emlem}{\end{mlem}}

\newcommand{\bscl}{\begin{scl}}
\newcommand{\escl}{\end{scl}}

\newcommand{\bcons}{\begin{conjs}}
\newcommand{\econs}{\end{conjs}}

\newcommand{\bprop}{\begin{prop}}
\newcommand{\eprop}{\end{prop}}

\newcommand{\br}{\begin{rem}}
\newcommand{\er}{\end{rem}}
\newcommand{\brs}{\begin{rems}}
\newcommand{\ers}{\end{rems}}
\newcommand{\bo}{\begin{obser}}
\newcommand{\eo}{\end{obser}}
\newcommand{\bos}{\begin{obsers}}
\newcommand{\eos}{\end{obsers}}
\newcommand{\bpf}{\begin{pf}}
\newcommand{\epf}{\end{pf}}
\newcommand{\ba}{\begin{array}}
\newcommand{\ea}{\end{array}}
\newcommand{\beq}{\begin{eqnarray}}
\newcommand{\beqq}{\begin{eqnarray*}}
\newcommand{\eeq}{\end{eqnarray}}
\newcommand{\eeqq}{\end{eqnarray*}}

\newcommand{\ds}{\displaystyle}

\newcounter{minutes}
\divide\time by 60
\newcounter{hours}
\multiply\time by 60 \addtocounter{minutes}{-\time}

\begin{document}

\bibliographystyle{amsplain}
\title []
{Advancements in Log-P-Analytic Functions: Landau-Type Theorems and Their Refinements} 

\def\thefootnote{}
\footnotetext{ \texttt{\tiny File:~\jobname .tex,
          printed: \number\day-\number\month-\number\year,
          \thehours.\ifnum\theminutes<10{0}\fi\theminutes}
} \makeatletter\def\thefootnote{\@arabic\c@footnote}\makeatother

\author{Hanghang Zhao}
\address{H.H. Zhao, School of Mathematical Sciences, South China Normal University, Guangzhou, Guangdong 510631, China.}
\email{zhh18888106641@163.com}

\author{Ming-Sheng Liu ${}^{~\mathbf{*}}$}
\address{M.S. Liu, School of Mathematical Sciences, South China Normal University, Guangzhou, Guangdong 510631, China.}
\email{liumsh@scnu.edu.cn}

\author{Kit Ian Kou}
\address{K.I. Kou, Department of Mathematics, Faculty of Science and Technology, University of Macau, Taipa, Macao, China}
\email{kikou@um.edu.mo} 

\subjclass[2000]{30C99, 31A30} 
\keywords{Landau-type theorem; Bloch theorem; poly-analytic function; log-$p$-analytic function;  univalent.
  \\
${}^{\mathbf{*}}$ Correspondence should be addressed to Ming-Sheng Liu 
}


\begin{abstract}
This work begins by introducing the groundbreaking concept of log-p-analytic functions. Following this introduction, we proceed to delineate four distinct formulations of Landau-type theorems, specifically crafted for the domain of poly-analytic functions. Among these, two theorems are distinguished by their exactitude, and a third theorem offers a refinement to the existing work of Abdulhadi and Hajj. Concluding the paper, we present four specialized versions of Landau-type theorems applicable to a subset of bounded log-p-analytic functions, resulting in the derivation of two precise outcomes.
\end{abstract}

\maketitle
\pagestyle{myheadings}
\markboth{H.H. Zhao, M.S. Liu and K.I. Kou}{Landau-type theorems of log-$p$-analytic functions}

\section{Introduction and Preliminaries}\label{HLP-sec1}
In the domain of classical complex analysis, a significant unresolved question pertains to the precise determination of the Bloch constant for functions that are analytic within the unit disk. Chen et al., in their work cited as \cite{CGH2000}, have examined a parallel issue, that of assessing the Bloch constant for planar harmonic mappings. This line of inquiry, propelled by the seminal contributions from \cite{CGH2000}, has since been pursued by various researchers, resulting in notable enhancements to the pre-existing Landau-type theorems. A detailed exposition of these improvements is scheduled for a later segment of this manuscript. The present paper extends this discourse by examining poly-analytic and log-$p$-analytic functions, within which we introduce a suite of newly articulated, exact Landau-type theorems.

\subsection{Definitions and Notations}
Suppose that $p$ is a positive integer and $F$ is a  continuous complex-valued function defined in a domain $D\subset \mathbb{C}$. Then $F$ is said to be poly-analytic function of order $p$ if it satisfies the generalized Cauchy-Riemann equations
$$
\frac{\partial ^{p}F\left( z \right)}{\partial \bar{z}^{p}}=0,\, z\in D.
$$

When $p = 1$, the function $F$ is analytic in $D$. When $p = 2$, the function $F$ is called bianalytic  in $D$. It is well-known that every poly-analytic function $F$ in a simply connected domain $D$  has the representation (cf. \cite{AH2022})
\begin{eqnarray}
F=\sum_{k=0}^{p -1}{\bar{z}^k}A_k\left( z \right),
\label{liu11}
\end{eqnarray}
where all $A_k\left( z \right)$ are analytic functions  in $D$ for $k=0, \cdots , p -1$.

The concept of poly-analytic functions initially emerged from the work of the Russian mathematician Kolossov, who was engaged in research pertaining to the mathematical theory of elasticity. Since their inception, poly-analytic functions have been the subject of extensive scholarly inquiry. A significant milestone in this field was achieved by Agranovsky, who in \cite{A2011}, provided a novel characterization of these functions. Progressing further, in 2011, Li and Wang \cite{LW2012} introduced the concept of log-p-harmonic mappings, defining a mapping \( f \) as log-$p$-harmonic if the logarithm of \( f \) qualifies as a $p$-harmonic mapping. Drawing inspiration from the work presented in \cite{LW2012}, this paper introduces the innovative concept of log-$p$-analytic functions, delineating their properties and implications in the broader mathematical landscape.

A function $f$ is {\it log-$p$-analytic} in a domain $D\subset \mathbb{C}$ if $\log f$ is poly-analytic function of order $p$ in $D$. Here and in what follows, ``$\log$'' is taken to be its principal branch such that $\log 1=0$.

When $p = 1$, the function $f$ is called log-analytic. When $p = 2$, the function $f$ is called log-bianalytic. Then it follows from (\ref{liu11}) that $f$ is log-$p$-analytic in a simply connected domain $D$ if and only if $f$ can be written as
$$
f(z)=\prod\limits_{k=0}^{p-1}{\left( a_k\left( z \right) \right)}^{\bar{z}^k},
$$
where $a_k(z)$ is log-analytic function for $k\in \left\{ 0,1,\ldots, p-1 \right\}$.

As per the findings reported in \cite{B1997}, it is established that a poly-analytic complex function possesses poly-harmonic properties; however, the converse of this statement does not hold. Similarly, in the realm of log-$p$-analytic complex functions, it is observed that they exhibit log-$p$-harmonic characteristics, yet the reciprocal of this assertion is not universally valid.

For $r>0$, let $U_r=\left\{ z\in \mathbb{C}:\left| z \right|<r \right\}$ denote the disk with center at the origin and radius $r$, in particular, $U:=U_1$.
For a complex-valued function $f$ in the unit disk $U$, its Jacobian $J_f(z)$ is given by $J_f(z)=|f_z(z)|^2-|f_{\overline{z}}(z)|^2$.
We say that a harmonic mapping $f$ is locally univalent and sense-preserving if and only if its Jacobian $J_f(z)>0$ for $z\in D$ (cf. \cite{L1936}). 
For continuously differentiable function $f$, let
\begin{eqnarray*}
\Lambda_{f}(z)
=|f_{z}(z)|+|f_{\overline{z}}(z)| ~\mbox{ and }~ \lambda_{f}(z)=\big ||f_{z}(z)|-|f_{\overline{z}}(z)|\big |.
\end{eqnarray*}

\subsection{Landau and Bloch theorems}
The classical theorem of Landau states that if $f$ is analytic in the unit disk $U$ with $f(0)=0, f'(0)=1$ and $|f(z)|<M$ in $U$ for some $M>1$, then $f$ is univalent in $U_{r _0}$ and $f(U_{r _0})$ contains a disk $U_{R_0}$ with $r _0=1/(M+\sqrt{M^2-1}),\, R_0=Mr_0^2$.
This result is sharp, with the extremal function $f_0(z)=M z \frac{1-Mz}{M-z}$ (see \cite{L2009S}). The Bloch theorem asserts the existence of a positive constant number $b$ such that if $f$ is an analytic function on the unit disk $U$ with $f'(0) = 1$, then $f(U)$ contains a schlicht disk of radius $b$, that is, a disk of radius $b$ which is the univalent image of some region in $U$. The supremum of all such constants $b$ is called the Bloch constant (see \cite{{CGH2000}} ).

In the year 2000, Chen et al. \cite{CGH2000}, under appropriate conditions, initially derived three variations of Landau-type theorems for a specific class of bounded harmonic mappings within the unit disk. Subsequent to this seminal work, numerous scholars have delved into the exploration of Landau-type theorems for harmonic mappings, further refining the initial findings (refer to \cite{cg,CPR2014,cpw2011C,CPW2011B,DN2004,G2006,Huang2014,L2009S,LC2018}). Concurrently, the attention of the academic community has also been drawn to the Landau-type theorems in the context of biharmonic mappings, as evidenced by the contributions in \cite{CPW2009, cpw2011C, L2008C,LXY2017,ZL2013}. Notably, in recent years, there has been a marked advancement in the acquisition of precise results for normalized bounded harmonic mappings \cite{Huang2014,L2009S,LC2018}, as well as for normalized bounded biharmonic mappings \cite{LL2019}. In the year 2024, Liu and Ponnusamy \cite{LP2024} achieved a significant milestone by presenting two exact formulations of Landau-type theorems for bianalytic mappings, as detailed below.

{\bf Theorem A.}\,{\rm (\cite[Theorem 2.1]{LP2024})}\quad
  Suppose that $\Lambda_1\geq0$ and $\Lambda_2>1$. Let $F$ be a bi-analytic functions and $F(z)=\bar{z}G(z)+H(z)$, where $H(z)$ and $G(z)$ are analytic, $H(0)=G(0)=0$, $H'(0)=1$,
  $|G'(z)|\leq\Lambda_1$ and $|H'(z)|<\Lambda_2$ for all $z\in U$. Then $F$ is univalent in $U_{r_1}$, and $F\left( U_{r_1} \right) $ contains a schlicht disk $U_{R_1}$, where
  \begin{eqnarray}
  r_1=\frac{2\Lambda_2}{\Lambda_2(2\Lambda_1+\Lambda_2)+\sqrt{\Lambda_2^2(2\Lambda_1+\Lambda_2)^2-8\Lambda_1\Lambda_2}},
  \end{eqnarray}
  and
  \begin{eqnarray}
  R_1= \Lambda_2^2r_1+\left(\Lambda_2^3-\Lambda_2\right)\ln\left(1-\frac{r_1}{\Lambda_2}\right)-\Lambda_1 r_1^2.
  \end{eqnarray}

 This result is sharp.

{\bf Theorem B.}\,{\rm (\cite[Theorem 2.2]{LP2024})}\quad
  Suppose that $\Lambda\geq 0$. Let $F$ be a bi-analytic functions and $F(z)=\bar{z}G(z)+H(z)$, where $H(z)$ and $G(z)$ are analytic, $H(0)=G(0)=0$, $H'(0)=1$,
  $|G'(z)|\leq\Lambda$, and $|H(z)|< 1$ or $|H'(z)|\leq 1$ for all $z\in U$. Then is univalent in $U_{r_2}$, and $F\left( U_{r_2} \right) $ contains a schlicht disk $U_{R_2}$, where
  $$
  r_2=\left \{ \begin{array}{rl}
  1 &\mbox{ when $\ds 0\leq\Lambda\leq\frac{1}{2}$}\\
  \ds \frac{1}{2\Lambda} &\mbox{ when $\ds \Lambda>\frac{1}{2}$},
  \end{array}
  \right .
  $$
  and $R_2=r_2-\Lambda r_2^2$. This result is sharp too.

In the year 2022, Abdulhadi and Hajj introduced a non-sharp variant of the Landau-type theorem, specifically tailored for a subset of bounded poly-analytic functions. 

{\bf Theorem C.}\,{\rm (\cite[Theorem 1]{AH2022})} \quad
  Let $F$ be a ploy-analytic function of order $p$, where all $A_k\left( z \right)$ are analytic, satisfying $A_k\left( 0 \right) =0,A_{k}'\left( 0 \right) =1$ and $\left| A_{k}\left( z \right) \right|\le M, M>1$ for all $k$ and $z\in U$.
  Then $F\left( z \right) $ is univalent in $U_{r_3}$, and $F\left( U_{r_3} \right) $ contains a schlicht disk $U_{R_3}$ , where $r_3$ is the unique root in $\left( 0,1 \right)$ of the following equation:
  $$
  1-M\left( \frac{r\left( 2-r \right)}{\left( 1-r \right) ^2}+\sum_{k=1}^{p -1}{\frac{r^{k}\left( 1+k-kr \right)}{\left( 1-r \right) ^2}} \right) =0
  $$
  and
  $$
  R_3=r_3-r_3^{2}\left( \frac{1-r_3^{p -1}}{1-r_3} \right) -M\sum_{k=0}^{p -1}{\frac{r_3^{k+2}}{1-r_3}}.
  $$


Theorem C is not sharp, a couple of natural questions arise.

\bprob\label{HLP-prob1}
Can we establish the sharp versions of Landau-type theorems for poly-analytic functions by extending Theorems A and B?  Whether
we can improve Theorem C?
\eprob

\bprob\label{HLP-prob2}
Can we establish several analogous versions of Landau-type theorems for certain bounded log-$p$-analytic functions? Moreover, can we establish several sharp results?
\eprob

The structure of this paper is delineated as follows. Section \ref{HLP-sec2} encompasses the articulation of eight theorems, among which one theorem offers an enhancement to Theorem C. Specifically, this section presents four precise renditions of Landau-type theorems applicable to both poly-analytic and log-$p$-analytic functions. These outcomes affirmatively address the inquiries posed in Problems \ref{HLP-prob1} and \ref{HLP-prob2}. Section \ref{HLP-sec3} introduces a pair of lemmas essential for substantiating the principal findings detailed in Section \ref{HLP-sec4}. The conclusion is given in Section \ref{HLP-sec5}.

\section{Statement of Main Results and Remarks}\label{HLP-sec2}

We initially articulate a definitive version of the Landau-type theorem, which is specifically crafted for a particular subclass of bounded poly-analytic functions. 
\bthm\label{HLP-th1}
Suppose that $\Lambda_0>1$ and $\varLambda _1,\cdots ,\varLambda _{p -1}\ge 0$. Let
$$F(z)=\sum_{k=0}^{p -1}{\bar{z}^k}A_k\left( z \right) $$
be a ploy-analytic function of order $p$, where all $A_k\left( z \right)$ are analytic, satisfying $A_k\left( 0 \right) =0,\, A_{0}'\left( 0 \right) =1$ and $|A'_0(z)|<\varLambda _0$, $\left| A_{k}'\left( z \right) \right|\le \varLambda _k, k\in \left\{ 1,\cdots , p -1 \right\} $ for all $z\in U$.
Then $F\left( z \right) $ is univalent in $U_{\rho_1}$ and $F\left( U_{\rho_1} \right) $ contains a schlicht disk $U_{\sigma_1}$, where $\rho_1$ is the unique root in $\left( 0,1 \right)$ of the following equation:
\begin{eqnarray}
  \varLambda _0\frac{1-\varLambda _0 r}{\varLambda _0-r}-\sum_{k=1}^{p -1}{\left( k+1 \right)}\varLambda _kr^k=0,
  \label{liu21}
\end{eqnarray}
and
\begin{eqnarray}
\sigma_1=\Lambda_0^2 \rho_1-\sum_{k=1}^{p -1}{\varLambda _k}\rho_1^{k+1}+\left(\Lambda_0^3-\Lambda_0\right)\ln\bigg(1-\frac{\rho_1}{\Lambda_0}\bigg).
\label{liu22}
\end{eqnarray}
This result is sharp, with an extremal function given by
\begin{align}
  F_1(z) &=\varLambda _0\int_{\left[ 0, z \right]}{\frac{\frac{1}{\varLambda _0}-z}{1-\frac{z}{\varLambda _0}}}dz-\sum_{k=0}^{p -1}{\varLambda _k}\bar{z}^k z \nonumber\\
   &=\Lambda_0^2 z-\sum_{k=1}^{p -1}{\varLambda _k}\bar{z}^kz+\left(\Lambda_0^3-\Lambda_0\right)\ln\bigg(1-\frac{z}{\Lambda_0}\bigg).
   \label{liu23}
  \end{align}

\ethm

Next, for the case $\Lambda_0=1$ and $\varLambda _1,\cdots ,\varLambda _{p -1}\ge 0$,  we will prove the following sharp version of Landau-type theorem for certain subclass of bounded poly-analytic functions.

\bthm \label{HLP-th2}
Suppose that $\varLambda _1,\cdots ,\varLambda _{p -1}\ge 0$. Let
$$F(z)=\sum_{k=0}^{p -1}{\bar{z}^k}A_k\left( z \right) $$
be a ploy-analytic function of order $p$, where all $A_k\left( z \right)$ are analytic, satisfying $A_k\left( 0 \right) =0, A_{0}'\left( 0 \right) =1$ $\left| A_{k}'\left( z \right) \right|\le \varLambda _k\ , k\in \left\{ 1,2,\cdots ,p -1 \right\} $ and $|A_0(z)|< 1$ or $|A_{0}'(z)|\leq 1$ for all $z\in U$.
Then $F\left( z \right) $ is univalent in $U_{\rho_2}$ and $F\left( U_{\rho_2} \right) $ contains a schlicht disk $U_{\sigma_2}$, where
\begin{eqnarray}
\rho _2=\left\{ \begin{array}{l}
	1,\ \mbox{ if }\sum\limits_{k=1}^{p -1}{\left( k+1 \right) \varLambda _k}\le 1 ,\\
	\rho _{2}',\ \mbox{ if }\sum\limits_{k=1}^{p -1}{\left( k+1 \right) \varLambda _k}>1 ,\\
\end{array} \right.
\label{liu24}
\end{eqnarray}
and $\rho _{2}'$ is the unique root in $\left( 0,1 \right)$ of the following equation:
\begin{eqnarray}
  1-\sum_{k=1}^{p -1}{\left( k+1 \right)}\varLambda _kr^k=0,
\end{eqnarray}
and
\begin{eqnarray}
  \sigma_2= \rho_2-\sum_{k=1}^{p -1}{\varLambda _k}\rho_2^{k+1}.
  \label{liu26}
  \end{eqnarray}
This result is sharp, with an extremal function given by
\begin{eqnarray}
  F_2(z) = z-\sum_{k=1}^{p -1}{\varLambda _k}\bar{z}^kz.
  \label{liu27}
   \end{eqnarray}
\ethm

\br\label{HLP-re1}
From \cite{B1997}, we know that a poly-analytic function is polyharmonic, hence, we conclude the following corollaries from Theorems \ref{HLP-th1} and \ref{HLP-th2}.
\er

\bcor\label{HLP-cor1}
Suppose that $\Lambda_0>1$ and $\varLambda _1,\cdots ,\varLambda _{p -1}\ge 0$. Let $F=\sum\limits_{k=0}^{p -1}{\bar{z}^k}A_k\left( z \right) $ be a polyharmonic mappings, where all $A_k\left( z \right)$ are analytic in $U$ 
with $A_k\left( 0 \right) =0, A_{0}'\left( 0 \right) =1$.
\begin{enumerate}
\item If $\left| A_{k}'\left( z \right) \right|\le \varLambda _k, k\in \left\{ 1,2,\cdots ,p -1 \right\}$ and $|A'_0(z)|<\varLambda _0$ for all $z\in U$, then $F\left( z \right) $ is univalent in $U_{\rho_1}$, and $F\left( U_{\rho_1} \right) $ contains a schlicht disk $U_{\sigma_1}$,
where $\rho_1$  and $\sigma_1$ are given by \eqref{liu21} and \eqref{liu22}, respectively. This result is sharp, with an extremal function  $F_1(z)$ given by \eqref{liu23}.

\item If $\left| A_{k}'\left( z \right) \right|\le \varLambda _k\ ,k\in \left\{ 1,2,\cdots ,p -1 \right\} $ and $|A_0(z)|< 1$ or $|A_{0}'(z)|\leq 1$ for all $z\in U$, then $F\left( z \right) $ is univalent in $U_{\rho_2}$, and $F\left( U_{\rho_2} \right) $ contains a schlicht disk $U_{\sigma_2}$, where $\rho_2$  and $\sigma_2$ are given by \eqref{liu24} and \eqref{liu26}, respectively. This result is sharp, with an extremal function  $F_2(z)$ given by \eqref{liu27}.

\end{enumerate}
\ecor

Now we improve Theorem C by establishing the following theorem.

\bthm\label{HLP-th3}
Suppose that $M _0,\cdots ,M _{p -1}> 1$. Let
$$F(z)=\sum_{k=0}^{p -1}{\bar{z}^k}A_k\left( z \right) $$be a ploy-analytic function of order $p$,
where all $A_k\left( z \right)$ are analytic, satisfying $A_k\left( 0 \right) =0,A_{k}'\left( 0 \right) =1$ and $\left| A_{k}\left( z \right) \right|\le M_k,  k\in \left\{ 0,1,\cdots ,p -1 \right\} $ for all $z\in U$.
Then $F\left( z \right) $ is univalent in $U_{\rho_3}$, and $F\left( U_{\rho_3} \right) $ contains a schlicht disk $U_{\sigma_3}$ , where $\rho_3$ is the unique root in $\left( 0,1 \right)$ of the following equation:
\begin{eqnarray}
1-\sum_{k=0}^{p -1}{\left( M_k-\frac{1}{M_k} \right)}\frac{r^{k+1}\left( 2-r +k(1-r) \right)}{\left( 1-r \right) ^2}-\sum_{k=1}^{p -1}{\left( k+1 \right)}r^k=0,
\label{liu28}
\end{eqnarray}
and
\begin{eqnarray}
\sigma_3=\rho_3-\rho _{3}^{2}\left( \frac{1-\rho _{3}^{p -1}}{1-\rho _3} \right) -\sum_{k=0}^{p -1}{\left( M_k-\frac{1}{M_k} \right)}\frac{\rho _{3}^{k+2}}{1-\rho _3}.
\label{liu29}
\end{eqnarray}
\ethm

\br
If we set $M_k=1$ for all $k$ in Theorem \ref{HLP-th3}, then it is clear that $A_k(z)=z$ $k\in \left\{ 0,1,\cdots p -1 \right\} $ by Schwarz lemma. Thus, $\rho_3$ is the unique root of $\sum\limits_{k=1}^{p -1}{\left( k+1 \right)}r^k=1$ and $\sigma_3=\rho_3-\rho _{3}^{2}\left(\displaystyle \frac{1-\rho _{3}^{p -1}}{1-\rho _3} \right)$ are sharp, with an extremal function  $F_3(z)=z+|z|^2\sum\limits_{k=0}^{p -2}{\bar{z}^k}$. Moreover, if we set $M_k = M$ for all $k$ in Theorem \ref{HLP-th3}, then one can easily gets an improved version of Theorem C.
\er

If we replace the condition $``| A_0\left( z \right) | \le M_0$  for all $z\in U"$ by the conditions $ ``| A_{0}'\left( z \right) |<\varLambda $ for all $z\in U"$ in Theorem  \ref{HLP-th3}, then, we conclude the following theorem.

\bthm\label{HLP-th4}
Suppose that $\Lambda>1$ and $M _1,\cdots ,M _{p -1}> 1$. Let
$$F(z)=\sum_{k=0}^{p -1}{\bar{z}^k}A_k\left( z \right) $$
be a ploy-analytic function of order $p$, where all $A_k\left( z \right)$ are analytic, satisfying $A_k\left( 0 \right) =0, A_{k}'\left( 0 \right) =1$ and $| A_{0}'\left( z \right) |< \varLambda$, $\left| A_{k}\left( z \right) \right|\le M_k\ ,k\in \left\{ 1,2,\cdots ,p -1 \right\} $ for all $z\in U$.
Then $F\left( z \right) $ is univalent in $U_{\rho_4}$ and $F\left( U_{\rho_4} \right) $ contains a schlicht disk $U_{\sigma_4}$, where $\rho_4$ is the unique root in $\left( 0,1 \right)$ of the following equation:
\begin{eqnarray}
\quad\varLambda \frac{1-\varLambda r }{\varLambda -r}-\sum_{k=1}^{p -1}{\left( M_k-\frac{1}{M_k} \right)}\frac{r^{k+1}\left( 2-r+k(1-r) \right)}{\left( 1-r \right) ^2}-\sum_{k=1}^{p -1}{\left( k+1 \right)}r^k=0,
\label{liu210}
\end{eqnarray}
and
\begin{eqnarray}
 \qquad \quad\sigma_4=\Lambda^2 \rho _{4}+(\Lambda^3-\Lambda)\ln\left(1-\frac{\rho _{4}}{\Lambda}\right)-\rho _{4}^{2}\left( \frac{1-\rho _{4}^{p -1}}{1-\rho _4} \right) -\sum_{k=1}^{p -1}{\left( M_k-\frac{1}{M_k} \right)}\frac{\rho _{4}^{k+2}}{1-\rho _4}.
\label{liu211}
\end{eqnarray}
\ethm

\br
Theorem 2.6 in \cite{LP2024} is just a special case of Theorem \ref{HLP-th3} when $p = 2$, and letting $p = 2$ in Theorem \ref{HLP-th4}, we obtain a new version of Landau-type theorem of bi-analytic in $U$ as follows:
\er

\bcor\label{HLP-cor2}
Let $F(z)=\bar{z}G(z)+H(z)$ be a bi-analytic function of the unit disk $U$, where $G(z)$ and $H(z)$ are analytic in $U$, satisfying $G(0)=H(0)=0$, $G'(0)=H'(0)=1, |G(z)|\leq M$ and $|H'(z)|\leq \varLambda$ for all $z\in U$. Then $F(z)$ is univalent in the disk $U_{\rho_4'}$ and $F(U_{\rho_4'})$ contains a schlicht disk $U_{\sigma_4'}$, where $\rho_4'$ is the unique root in $(0,\,1)$ of the equation
\begin{eqnarray*}
  \varLambda \frac{1-\varLambda r }{\varLambda -r}-{\left( M-\frac{1}{M} \right)}\frac{r^2\left( 3-2r \right)}{\left( 1-r \right) ^2}-2r=0,
\end{eqnarray*}
and
\begin{eqnarray*}
  \sigma_4'=\Lambda^2 \rho _4'+(\Lambda^3-\Lambda)\ln\left(1-\frac{\rho _4'}{\Lambda}\right)-\rho _4'^{2} -{\left( M-\frac{1}{M} \right)}\frac{\rho _4'^{3}}{1-\rho _4'}.
\end{eqnarray*}
\ecor

In conclusion, we delineate the exact formulation of Landau-type theorems for a select category of log-$p$-analytic functions, leveraging Theorem \ref{HLP-th1} and Lemma \ref{HLP-lem3}. Notably, this approach culminates in the determination of the precise value \( r_1 = \sinh\sigma_1 \) for the Bloch constant associated with specific log-$p$-analytic mappings within the unit disk \( U \).

\bthm\label{HLP-th5}
Suppose that $\Lambda_0>1$ and $\varLambda _1,\cdots ,\varLambda _{p-1}\ge 0$. Let$$f(z)=\prod\limits_{k=0}^{p-1}{\left( a_k\left( z \right) \right)}^{\bar{z}^k}$$be a log-p-analytic function, satisfying $f\left( 0 \right) =\lambda_f\left( 0 \right)=1$.
Suppose that for each $k\in \left\{ 0,1,\cdots, p-1 \right\}$ we have that
\begin{itemize}
  \item [(i)] $a_k(z)$ is log-analytic in $U$, $a_k(0)=a_k'(0)=1$ and $A_k(z):=\log a_k(z)$;
  \item [(ii)] for each $k\in \left\{ 1,\cdots, p-1 \right\}$, $\left| A_{k}'\left( z \right) \right|\le \varLambda _k$, and $|A_{0}'(z)|< \varLambda _0$ for all $z \in U$.
\end{itemize}

Then $f\left( z \right) $ is univalent in $U_{\rho_1}$, where $\rho_1$ is the unique root in $\left( 0,1 \right)$ of the equation \eqref{liu21}. Moreover, the range $f\left( U_{\rho_1} \right) $ contains a schlicht disk $U\left( w_1,r_1 \right) $, where $\sigma_1$ is given by \eqref{liu22}, and
\begin{eqnarray}
  w_1=\cosh\sigma_1,\quad r_1=\sinh\sigma_1.
  \label{liu210}
\end{eqnarray}
This result is sharp.
\ethm

Utilizing Theorem \ref{HLP-th2} and Lemma \ref{HLP-lem3}, and employing a methodology analogous to that used in the proof of Theorem \ref{HLP-th5}, we derive the specific formulation of the Landau-type theorem for a subset of log-$p$-analytic mappings. Specifically, this derivation yields the exact value \( r_2 = \sinh\sigma_2\) of the Bloch constant for these particular log-$p$-analytic mappings within the unit disk \( U \).

\bthm\label{HLP-th6}
Suppose that  $\varLambda _1,\cdots ,\varLambda _{p-1}\ge 0$. Let
$$f(z)=\prod\limits_{k=0}^{p-1}{\left( a_k\left( z \right) \right)}^{\bar{z}^k}$$
be a log-p-analytic function, satisfying $f\left( 0 \right)=\lambda_f\left( 0 \right) =1$.
Suppose that for each $k\in \left\{ 0,1,\cdots, p-1 \right\}$ we have that
\begin{itemize}
  \item [(i)] $a_k(z)$ is log-analytic in $U$, $a_k(0)=a_k'(0)=1$ and $A_k(z):=\log a_k(z)$;
  \item [(ii)] for each $k\in \left\{ 1,\cdots, p-1 \right\}$, $\left| A_{k}'\left( z \right) \right|\le \varLambda _k$, and $|A_{0}'(z)|\le 1$ for all $z \in U$.
\end{itemize}

Then $f\left( z \right) $ is univalent in $U_{\rho_2}$, where $\rho_2$ is the unique root in $\left( 0,1 \right)$ of the equation\eqref{liu24}.
Moreover, the range $f\left( U_{\rho_2} \right) $ contains a schlicht disk $U\left( w_2,r_2 \right) $, where $\sigma_2$ is given by \eqref{liu26}, and
\begin{eqnarray}
  w_2=\cosh\sigma_2,\quad r_2=\sinh\sigma_2.
\end{eqnarray}
This result is sharp.
\ethm

Employing Theorems \ref{HLP-th3} and \ref{HLP-th4}, along with Lemma \ref{HLP-lem3}, and adhering to the methodological approach established in our verification of Theorem \ref{HLP-th5}, we present the ensuing pair of theorems.

{\bf Theorem 7.}\quad Suppose that $M_{0}^{*}, \cdots, M_{p-1}^{*}>1$. Let
$$f(z)=\prod\limits_{k=0}^{p-1}{\left( a_k\left( z \right) \right)}^{\bar{z}^k}$$
be a log-p-analytic function, satisfying $f\left( 0 \right)=\lambda_f\left( 0 \right) =1$. Suppose that for each $k\in \left\{ 0,1,\cdots, p-1 \right\}$ we have that

\begin{itemize}
  \item [(i)] $a_k(z)$ is log-analytic in $U$, $a_k(0)=a_k'(0)=1$, and $A_k(z):=\log a_k(z)$;
  \item [(ii)] $\left| a_{k}\left( z \right) \right|\le M_{k}^{*}$ for all $z \in U$.
\end{itemize}
Then $f\left( z \right) $ is univalent in $U_{\rho_3}$, where $M_k=\log M_{k}^{*}+\pi$ and $\rho_3$ is the unique root in $\left( 0,1 \right)$ of the equation \eqref{liu28}.
Moreover, the range $f\left( U_{\rho_3} \right) $ contains a schlicht disk $U\left( w_3,r_3 \right) $, where $\sigma_3$ is given by \eqref{liu29}, and
\begin{eqnarray}
  w_3=\cosh\sigma_3,\quad r_3=\sinh\sigma_3.
  \label{liu212}
\end{eqnarray}


{\bf Theorem 8.}\quad Suppose that $\Lambda>1$ and $M_{1}^{*},\cdots ,M_{p-1}^{*}>1$. Let
$$f(z)=\prod\limits_{k=0}^{p-1}{\left( a_k\left( z \right) \right)}^{\bar{z}^k}$$
be a log-p-analytic function, satisfying $f\left( 0 \right)=\lambda_f\left( 0 \right) =1$.
Suppose that for each $k\in \left\{ 0,1,\cdots, p-1 \right\}$ we have that
\begin{itemize}
  \item [(i)] $a_k(z)$ is log-analytic in $U$, $a_k(0)=a_k'(0)=1$, and $A_k(z):=\log a_k(z)$;
  \item [(ii)] for each $k\in \left\{ 1,\cdots, p-1 \right\}$, $\left| a_{k}\left( z \right) \right|\le M_{k}^{*}$, and $|A_0'(z)|< \Lambda $ for all $z \in U$.
\end{itemize}
Then $f\left( z \right) $ is univalent in $U_{\rho_4}$, where $M_k=\log M_{k}^{*}+\pi$ and $\rho_4$ is the unique root in $\left( 0,1 \right)$ of the equation \eqref{liu210}. Moreover, the range $f\left( U_{\rho_4} \right) $ contains a schlicht disk $U\left( w_4,r_4 \right) $, where $\sigma_4$ is given by \eqref{liu211}, and
\begin{eqnarray}
  w_4=\cosh\sigma_4,\quad r_4=\sinh\sigma_4.
  \label{liu213}
\end{eqnarray}


\section{Key lemmas}\label{HLP-sec3}
To substantiate our principal findings, we rely on the forthcoming lemmas, which are pivotal in the formulation of the results that are detailed in Section \ref{HLP-sec4}.

\blem\label{HLP-lem1}{\rm (\cite{LP2024}) }
 Let $H(z)$ be a analytic functions satisfying $|H'(0)|=1$ and $|H'(z)|<\Lambda$ for all $z\in U$ and for some $\Lambda>1$.
\begin{enumerate}
\item For all $z_1,z_2\in U_r\, (0<r<1, z_1\neq z_2)$, we have
$$
|H(z_1)-H(z_2)|=\bigg|\int_{\gamma }H'(z)\,dz\bigg|\geq\Lambda\, \frac{1-\Lambda r}{\Lambda-r}\, |z_1-z_2|,
$$
where $\gamma =[z_1,z_2]$ denotes the closed line segment joining $z_1$ and $z_2$.

\item  For $z'\in \partial U_r\, (0<r<1)$ with $w'=H(z')\in H(\partial U_r)$ and $|w'|=\min\left\{|w|:\,w\in H\left(\partial U_r\right)\right\}$, set $\gamma _0=H^{-1}(\Gamma _0)$ and
$\Gamma _0= [0,w'] $ denotes the closed line segment joining the origin and $w'$. Then we have
$$
|H(z')| \geq\Lambda\int_{0}^{r}\frac{\frac{1}{\Lambda}-t}{1-\frac{t}{\Lambda}}\,dt=\Lambda^2 r+(\Lambda^3-\Lambda)\ln\left(1-\frac{r}{\Lambda}\right).
$$
\end{enumerate}
\elem

\blem\label{HLP-lem2}{\rm (\cite{LP2024}) }
 Suppose that $f$ is a analytic functions satisfying $|H(0)|=0$, $|H'(0)|=1$ and $|H(z)|<M$ for some $M>0$ and $f(z)=\sum_{n=1}^{\infty}a_nz^n$, then $M\geq 1$  and
\begin{equation*}
|a_n|\leq M-\frac{1}{M} ~\mbox{ for  $n=2, 3, \ldots .$} 
\end{equation*}
These inequalities are sharp, with the extremal functions $f_n(z)$, where
\begin{eqnarray*}
f_1(z)=z,\quad f_n(z)=Mz\frac{1-M z^{n-1}}{M-z^{n-1}}=z-\Big(M-\frac{1}{M}\Big)z^n-\sum\limits_{k=3}^\infty\frac{M^2-1}{M^{k-1}}z^{(n-1)(k-1)+1}
\end{eqnarray*}
for $n=2,3,\ldots$.
\elem

By employing arguments analogous to those utilized in the proof referenced in \cite[Lemma 2.4]{LL2021}, the subsequent lemma can be established, thus rendering the provision of detailed proof unnecessary.

\blem\label{HLP-lem3}
Suppose that $p$ is a positive integer and $0<\sigma<1,\, 0<\rho\leq1$. Let $f(z)$ be a log-p-analytic function satisfying $f\left( 0 \right) =\lambda_f\left( 0 \right) =1$.
Suppose that $f(z)$ is univalent in $U_{\rho}$ and $F(U_{\rho})\supset U_{\sigma}$, where $F(z)=\log f(z)$. Then the range $f(U_{\rho})$ contains a schlicht disk $U( w_{0}, r_{0}) = \{ w\in \mathbb{C} ||w- w_{0}|< r_{0}\}$, where
\begin{eqnarray}
  w_0=\cosh\sigma,\quad r_0=\sinh\sigma.
\end{eqnarray}
Moreover, if $\rho$ is the biggest univalent radius of $f(z)$, then the radius $r_0=\sinh\sigma$ is sharp.
\elem

\blem\label{HLP-lem4}
Suppose that $\Lambda_0>1$ and $\varLambda _1,\cdots ,\varLambda _{p -1}\ge 0$. Let
$$g_1\left( r \right)=\Lambda_0 \frac{1-\Lambda_0 r}{\Lambda_0-r}-\sum_{k=1}^{p -1}{\left( k+1 \right)}\varLambda _kr^k.$$
Then, the function $g_1\left( r \right)$ is strictly decreasing on $\left[ 0,1 \right]$, and there exists a unique real number $\rho_1\in (0, 1/\varLambda _0]$ such that $g_1\left(\rho_1 \right)=0$.
\elem

\bpf
By the hypothesis, we have that $\Lambda_0>1$, and 
$$
g_1'\left( r \right)=\Lambda_0 \frac{1-\Lambda_0^2}{(\Lambda_0-r)^2}-\sum_{k=1}^{p -1}k{\left( k+1 \right)}\varLambda _kr^{k-1}<0,\quad r\in [0, 1],
$$
this implies that $g_1\left( r \right)$ is strictly decreasing on $\left[ 0,1 \right]$.

As the function $g_1\left( r \right)$ is continuous on $\left[ 0, 1 \right] $, $g_1\left( 0 \right)=1>0$, and
$$
g_1\left( \frac{1}{\varLambda _0} \right) =-\sum_{k=1}^{p -1}{\left( k+1 \right)}\varLambda _k(\frac{1}{\varLambda _0} )^k\le 0,
$$
therefore, by the mean value theorem, we have that there is a unique real number $\rho_1$ in $(0,1/\varLambda _0]$ such that $g_1\left(\rho_1 \right)=0$.
\epf

\blem\label{HLP-lem5}
Suppose that $M _0,\cdots ,M _{p -1}> 1$. Let
$$g_2\left( r \right)=1-\sum_{k=0}^{p -1}{\left( M_k-\frac{1}{M_k} \right)}\frac{r^{k+1}\left( 2-r +k(1-r) \right)}{\left( 1-r \right) ^2}-\sum_{k=1}^{p -1}{\left( k+1 \right)}r^k.$$
Then, the function $g_2\left( r \right)$ is strictly decreasing on $\left[ 0,1 \right]$, and there exists a unique real number $\rho_3\in (0, 1)$ such that $g_2\left(\rho_3 \right)=0$.
\elem

\bpf
By the hypothesis, we have that $M _0,\cdots ,M _{p -1}> 1$, and 
\begin{eqnarray*}
g_2'\left( r \right)&=&-\sum_{k=0}^{p -1}{\left( M_k-\frac{1}{M_k} \right)}\frac{r^{k}\left(k(k+1)r^2-2k(k+2)r+(k+1)(k+2)\right)}{\left( 1-r \right)^3} \\ &&-\sum_{k=1}^{p -1}k{\left( k+1 \right)}r^{k-1}<0,\quad r\in [0, 1],
\end{eqnarray*}
this implies that $g_2\left( r \right)$ is strictly decreasing on $\left[ 0,1 \right]$.

As the function $g_2\left( r \right)$ is continuous on $\left[ 0, 1 \right] $, $g_2\left( 0 \right)=1>0$, and
$$
\lim\limits_{r\to 1^-}g_2\left(r\right) =-\infty,
$$
therefore, by the mean value theorem, we have that there is a unique real number $\rho_3$ in $(0,1)$ such that $g_2\left(\rho_3 \right)=0$.
\epf

Adopting the methodological approach parallel to that applied in the proofs of Lemmas \ref{HLP-lem4} and \ref{HLP-lem5}, the ensuing lemma can be validated.
\blem\label{HLP-lem6}
Suppose that $\Lambda>1$ and $M _1,\cdots ,M _{p -1}> 1$. Let
$$g_3\left( r \right)=\varLambda \frac{1-\varLambda r }{\varLambda -r}-\sum_{k=0}^{p -1}{\left( M_k-\frac{1}{M_k} \right)}\frac{r^{k+1}\left( 2-r +k(1-r) \right)}{\left( 1-r \right) ^2}-\sum_{k=1}^{p -1}{\left( k+1 \right)}r^k.$$
Then, the function $g_3\left( r \right)$ is strictly decreasing on $\left[ 0, 1 \right]$, and there exists a unique real number $\rho_4\in (0, 1)$ such that $g_3\left(\rho_4 \right)=0$.
\elem

\section{Proofs of the main results}\label{HLP-sec4}

\subsection{Proof of Theorem \ref{HLP-th1}}
By the hypothesis of Theorem \ref{HLP-th1}, we have that for each $k\in \left\{ 0,1,\cdots p -1 \right\} $, $A_k\left( z \right)$ is analytic in $U$ and $A_k\left( 0 \right) =0$, which implies that
\begin{eqnarray}
|A_k(z)|=\bigg|\int_{[0,z]}A_k'(z)\, dz\bigg|\leq \int_{[0,z]}|A_k'(z)|\,|dz|\leq \Lambda_k |z|.
\label{liu41}
\end{eqnarray}

We first prove that $F$ is univalent in the disk $U_{\rho_1}$. In fact, for any $z_1,z_2\in U_r\, (0<r<\rho_1$, $z_1\neq z_2)$, we obtain from Lemmas  \ref{HLP-lem1} and \ref{HLP-lem4} that
\begin{eqnarray*}
|F(z_2)-F(z_1)|&=&\left| \int_{\left[ z_1,z_2 \right]}{F_z\left( z \right)}dz+F_{\bar{z}}\left( z \right) d\bar{z} \right|\nonumber\\
&=&\left| \int_{\left[ z_1,z_2 \right]}{\left( \sum_{k=0}^{p -1}{\bar{z}^k}A_{k}'\left( z \right) dz+\sum_{k=1}^{p -1}{k\bar{z}^{k-1}}A_k\left( z \right) d\bar{z} \right)} \right| \nonumber\\
&\geq&\left| \int_{\left[ z_1,z_2 \right]}{A_{0}'\left( z \right) dz} \right|-\left| \int_{\left[ z_1,z_2 \right]}{\left( \sum_{k=1}^{p -1}{\bar{z}^k}A_{k}'\left( z \right) dz+\sum_{k=1}^{p -1}{k\bar{z}^{k-1}}A_k\left( z \right) d\bar{z} \right)} \right|\\
&\geq&|z_1-z_2|\left(\Lambda_0 \frac{1-\Lambda_0 r}{\Lambda_0-r}-\sum_{k=1}^{p -1}{\left( k+1 \right)}\varLambda _kr^k\right)\\
&>&|z_1-z_2|\left(\Lambda_0 \frac{1-\Lambda_0 \rho_1}{\Lambda_0-\rho_1}-\sum_{k=1}^{p -1}{\left( k+1 \right)}\varLambda _k\rho_1^k\right)=0.
\end{eqnarray*}
Thus $F(z_2)\ne F(z_1)$, which implies the univalency of $F$ in the disk $U_{\rho_1}$.

Next, we prove that $F(U_{\rho_1}) \supseteq U_{\sigma_1}$. We note that $F(0)=0$, for $z'\in \partial U_{\rho_1}$ with $w'=F(z')\in F(\partial U_{\rho_1})$ and $|w'|=\min\left\{|w|:\,w\in F\left(\partial U_{\rho_1}\right)\right\}$, by \eqref{liu41} and Lemma \ref{HLP-lem1}, we have that
\begin{eqnarray*}
|w'|&=&\left|\sum_{k=0}^{p -1}{\bar{z'}^k}A_k\left( z' \right)\right|\geq\left| A_0\left( z' \right) \right|-\sum_{k=1}^{p -1}{\rho _{1}^{k}\left| A_k\left( z' \right) \right|}\\
&\geq&\Lambda_0^2 \rho_1-\sum_{k=1}^{p -1}{\varLambda _k}\rho_1^{k+1}+\left(\Lambda_0^3-\Lambda_0\right)\ln\bigg(1-\frac{\rho_1}{\Lambda_0}\bigg)=\sigma_1,
\end{eqnarray*}
which implies that $F(U_{\rho_1})\supseteq U_{\sigma_1}$.

Now, we prove the sharpness of $\rho_1$ and $\sigma_1$. To this end, we consider the poly-analytic function $F_1(z)$ which is given by (\ref{liu23}). It is easy to verify that $F_1(z)$ satisfies the hypothesis of Theorem \ref{HLP-th1}, and thus, we have that $F_1(z)$ is univalent in $U_{\rho_1}$, and $F_1(U_{\rho_1}) \supseteq U_{\sigma_1}$.

To show that the radius $\rho_1$ is sharp, we need to prove that $F_1(z)$ is not univalent in $U_r$ for each $r\in (\rho_1, 1]$, Indeed let
\begin{eqnarray}
  g(x)=\Lambda_0^2 x-\sum_{k=1}^{p -1}{\varLambda _k}x^{k+1}+\left(\Lambda_0^3-\Lambda_0\right)\ln\bigg(1-\frac{x}{\Lambda_0}\bigg), \quad  x\in [0, 1].
  \label{liu42}
  \end{eqnarray}

Because the continuous function
\begin{eqnarray}
  g'(x)=\Lambda_0^2 -\sum_{k=1}^{p -1}(k+1){\varLambda _k}x^k-\left(\Lambda_0^3-\Lambda_0\right)(\frac{1}{\Lambda_0-x})
  \label{liu43}
  \end{eqnarray}
is strictly decreasing on $[0, 1]$, and
\begin{eqnarray*}
 g'(\rho_1)&=&\Lambda_0^2 -\sum_{k=1}^{p -1}(k+1){\varLambda _k}\rho_1^k-\left(\Lambda_0^3-\Lambda_0\right)(\frac{1}{\Lambda_0-\rho_1})\\
 &=&\Lambda_0^2 -\Lambda_0 \frac{1-\Lambda_0 r}{\Lambda_0-r}-\left(\Lambda_0^3-\Lambda_0\right)(\frac{1}{\Lambda_0-\rho_1})=0,
\end{eqnarray*}
we obtain that $g'(x)=0$ for $x\in [0, 1]$ if and only if $x=\rho_1$. So $g(x)$ is strictly increasing on $[0, \rho_1)$
and strictly decreasing on $[\rho_1, 1]$. Since $g(0)=0$, there is a unique real $\rho_1'\in (\rho_1, 1]$ such that $g(\rho_1')=0$ if $g(1)\leq 0$, and
\begin{eqnarray}
  \sigma_1=\Lambda_0^2 \rho_1-\sum_{k=1}^{p -1}{\varLambda _k}\rho_1^{k+1}+\left(\Lambda_0^3-\Lambda_0\right)\ln\bigg(1-\frac{\rho_1}{\Lambda_0}\bigg)=g(\rho_1)>g(0)=0.
\label{liu44}
\end{eqnarray}

For every fixed $r\in (\rho_1, 1]$, set $x_1=\rho_1+\varepsilon$, where
\begin{equation*}
\varepsilon=\left\{
\begin{array}{lll}
\ds\min\left\{\frac{r-\rho_1}{2}, \frac{\rho_1'-\rho_1}{2}\right\} && \mbox{ if } g(1)\leq 0,\\
\ds \frac{r-\rho_1}{2} && \mbox{ if } g(1)>0.
\end{array}
\right.
\end{equation*}

By the mean value theorem, there is a unique $\delta\in (0, \rho_1)$ such that $x_2:=\rho_1-\delta\in (0, \rho_1)$ and $g(x_1)=g(x_2)$.

Let $z_1=x_1$ and $z_2=x_2$. Then $z_1,\ z_2\in U_r$ with $z_1\neq z_2$ and observe that
\begin{eqnarray*}
F_1(z_1)=F_1(x_1)=g(x_1)=g(x_2)=F_1(z_2).
\end{eqnarray*}
Hence $F_1$ is not univalent in the disk $U_r$ for each $r\in (\rho_1, 1]$, and thus, the radius $\rho_1$ is sharp.

Finally, note that $F_1(0)=0$ and picking up $z'=\rho_1\in \partial U_{\rho_1}$, by (\ref{liu22}), (\ref{liu42}) and (\ref{liu44}), we have
\begin{eqnarray*}
|F_1(z')-F_1(0)|=|F_1(\rho_1)|=|g(\rho_1)|=g(\rho_1)=\sigma_1.
\end{eqnarray*}
Hence, the radius $\sigma_1$ of the schlicht disk is also sharp.\hfill $\Box$

\subsection{Proof of Theorem \ref{HLP-th2}}
By the hypothesis of Theorem \ref{HLP-th2}, we have that $A_0\left( z \right)$ is analytic, $A_0\left( 0 \right) =0$, $A'_0\left( 0 \right) =1$ and $|A_0(z)|< 1$ or $|A_{0}'(z)|\leq 1$ for all $z\in U$. By Schwarz's lemma, we get that $A_0(z)\equiv z$ in $U$. Thus, $F$ reduces to the form
$$
F=z+\sum_{k=1}^{p -1}{\bar{z}^k}A_k\left( z \right).
$$

Now we prove $F$ is univalent in the disk $U_{\rho_2}$. To this end, for any $z_1,z_2\in U_r\, (0<r<\rho_2)$ with $z_1\neq z_2$. By the hypothesis of this theorem, it follows from (\ref{liu41}) that $|A_k(z)|\leq \Lambda_k |z|$ in $U$ for all $k\in \left\{ 1, 2, \cdots, p -1 \right\}$. Consequently,
\begin{eqnarray*}
|F(z_1)-F(z_2)|&\geq&\left| \int_{\left[ z_1,z_2 \right]}{dz} \right|-\left| \int_{\left[ z_1,z_2 \right]}{\left( \sum_{k=1}^{p -1}{\bar{z}^k}A_{k}'\left( z \right) dz+\sum_{k=1}^{p -1}{k\bar{z}^{k-1}}A_k\left( z \right) d\bar{z} \right)} \right|\\
&\geq&|z_1-z_2|\left(1-\sum_{k=1}^{p -1}{\left( k+1 \right)}\varLambda _kr^k\right)>0,\nonumber
\end{eqnarray*}
which proves the univalency of $F$ in the disk $U_{\rho_2}$, where $\rho_2$ is given in the statement of the theorem.

Noticing that $A_k(0)=0$, for any $z=\rho_2 e^{i\theta}\in \partial U_{\rho_2}$, we have
\begin{eqnarray*}
|F(z)|\geq|z|-\sum_{k=1}^{p -1}{\rho _{2}^{k}\left| A_k\left( z' \right) \right|}\geq\rho_2-\sum_{k=1}^{p -1}{\varLambda _k}\rho_2^{k+1}=\sigma_2.
\end{eqnarray*}
Hence, $F(U_{\rho_2})$ contains a schlicht disk $U_{\sigma_2}$.

Finally, for $ F_2(z) = z-\sum\limits_{k=1}^{p -1}{\varLambda _k}\overline{z}^kz$, a direct computation verifies that $\rho_2$ and $\sigma_2$ are sharp. This completes the proof.\hfill $\Box$

\subsection{Proof of Theorem \ref{HLP-th3}}
By the hypothesis of Theorem \ref{HLP-th3}, we have that all $A_k\left( z \right)$ are analytic, $A_k\left( 0 \right) =0, A_{k}'\left( 0 \right) =1$ and $\left| A_{k}\left( z \right) \right|\le M_k$ 
for all $z\in U$. Then for $k\in \left\{ 0,1,\cdots ,p -1 \right\} $, we may write
$$
A_k\left( z \right) =\sum\limits_{n=1}^{\infty}{a_{k,n}}z^n,
$$
where $a_{0,1}=a_{1,1}=\cdots a_{p -1,1}=1$, and it follows from Lemma \ref{HLP-lem2} that
\begin{eqnarray}
|a_{k,n}|\leq M_k-\frac{1}{M_k} ~\mbox{ for all $n\geq 2$  and }\, k\in \left\{ 0,1,\cdots ,p -1 \right\}.
\label{liu45}
\end{eqnarray}

We first prove that $F$ is univalent in the disk $U_{\rho_3}$. Indeed, for all $z_1,z_2\in U_r\, (0<r<\rho_3$, $z_1\neq z_2)$, we obtain from (\ref{liu45}) and Lemma  \ref{HLP-lem5} that
\begin{align*}
&|F(z_{1})-F(z_{2})|=\left|\int_{\left[ z_1,z_2 \right]}{F_z\left( z \right)}dz+F_{\bar{z}}\left( z \right) d\bar{z}\right|\\
&=\left|\int_{\left[ z_1,z_2 \right]}(F_{z}(0)dz+F_{\bar{z}}(0)d\bar{z})+\int_{\left[ z_1,z_2 \right]}\left(  F_{z}(z)-F_{z}(0)\right)  dz+\left(  F_{\bar{z}}(z)-F_{\bar{z}}(0)\right)d\bar{z}\right| \\
&\geq |z_{2}-z_{1}|-\int_{\left[ z_1,z_2 \right]}{\sum\limits_{n=2}^{\infty}}n|a_{0,n}||z|^{n-1}|dz|-
\int_{\left[ z_1,z_2 \right]}{\sum_{k=1}^{p-1}}|\bar{z}^{k-1}| (|z||A_{k}'(z)|+k|A_{k}(z)|)|dz|\\
\end{align*}\begin{align*}
&\geq|z_{2}-z_{1}|\left(1-{\sum_{n=2}^{\infty}}n|a_{0,n}|r^{n-1}-{\sum_{k=1}^{p-1}}r^{k}{\sum_{n=1}^{\infty}}(n|a_{k,n}|+k|a_{k,n}|)r^{n-1}\right) \\
&\geq|z_{2}-z_{1}|\left(1-{\sum_{k=0}^{p-1}}r^{k}\left( M_k-\frac{1}{M_k} \right){\sum_{n=2}^{\infty}}(n+k)r^{n-1}-{\sum_{k=1}^{p-1}}(k+1)r^{k}\right) \\
&=|z_{2}-z_{1}|\left(1-{\sum_{k=0}^{p-1}}\left( M_k-\frac{1}{M_k} \right)\frac{r^{k+1}(2-r+k(1-r))}{(1-r)^{2}}-{\sum_{k=1}^{p-1}}(k+1)r^{k}\right)>0.
\end{align*}
This implies $F(z_1)\neq F(z_2)$, which proves the univalency of $F$ in the disk $U_{\rho_3}$.

Next, we prove that $F(U_{\rho_3}) \supseteq U_{\sigma_3}$. Indeed, note that $F(0)=0$ and for any $z'\in \partial U_{\rho_3}$ with $w'=F(z')\in F(\partial U_{\rho_3})$, it follows from (\ref{liu45}) that
\begin{eqnarray*}
|w'|&=&\left|\sum_{k=0}^{p -1}{\bar{z'}^k}A_k\left( z' \right)\right|=\left| a_{0,1}z'+\sum_{k=1}^{p -1}{a_{k,1}}\left| z' \right|^2\bar{z'}^{k-1}+\sum_{k=0}^{p -1}{\bar{z}^k\sum_{n=2}^{\infty}{a_{k,n}z'^n}} \right|\\
&\geq&\rho _3-\sum_{k=1}^{p -1}{\rho _{3}^{k+1}}-\sum_{k=0}^{p -1}{\left( M_k-\frac{1}{M_k} \right)}\sum_{n=2}^{\infty}{\rho _{3}^{n+k}}\\
&=&\rho_3-\rho _{3}^{2}\left( \frac{1-\rho _{3}^{p -1}}{1-\rho _3} \right) -\sum_{k=0}^{p -1}{\left( M_k-\frac{1}{M_k} \right)}\frac{\rho _{3}^{k+2}}{1-\rho _3}=\sigma_3,
\end{eqnarray*}
which implies that $F(U_{\rho_3})\supseteq U_{\sigma_3}$. 
\hfill $\Box$

\subsection{Proof of Theorem \ref{HLP-th4}}
We first prove that $F$ is univalent in the disk $U_{\rho_4}$. Indeed, for all $z_1,z_2\in U_r\, (0<r<\rho_4$, $z_1\neq z_2)$, it follows from Lemmas \ref{HLP-lem1} and \ref{HLP-lem6}, and (\ref{liu45}) that
\begin{align*}
  &|F(z_{1})-F(z_{2})|=\left|\int_{\left[ z_1,z_2 \right]}{F_z\left( z \right)}dz+F_{\bar{z}}\left( z \right) d\bar{z}\right|\\
  &=\left|\int_{\left[ z_1,z_2 \right]}A_{0}'(z)dz\right|-\left|\int_{\left[ z_1,z_2 \right]}{\left( \sum_{k=1}^{p -1}{\bar{z}^k}A_{k}'\left( z \right) dz+\sum_{k=1}^{p -1}{k\bar{z}^{k-1}}A_k\left( z \right) d\bar{z} \right)}\right|\\
  &\geq|z_{2}-z_{1}|\left(\varLambda \frac{1-\varLambda r }{\varLambda -r}-{\sum_{k=0}^{p-1}}r^{k}\left( M_k-\frac{1}{M_k} \right){\sum_{n=2}^{\infty}}(n+k)r^{n-1}-{\sum_{k=1}^{p-1}}(k+1)r^{k}\right) \\
   &=|z_{2}-z_{1}|\left(\varLambda \frac{1-\varLambda r }{\varLambda -r}-{\sum_{k=0}^{p-1}}\left( M_k-\frac{1}{M_k} \right)\frac{r^{k+1}(2-r+k(1-r))}{(1-r)^{2}}-{\sum_{k=1}^{p-1}}(k+1)r^{k}\right)>0.
  \end{align*}
  This implies $F(z_1)\neq F(z_2)$, which proves the univalency of $F$ in the disk $U_{\rho_4}$.

  Next, we prove that $F(U_{\rho_4}) \supseteq U_{\sigma_4}$. Indeed, note that $F(0)=0$ and for any $z'\in \partial U_{\rho_4}$ with $w'=F(z')\in F(\partial U_{\rho_4})$,   it follows from Lemma \ref{HLP-lem1} and (\ref{liu45}) that
  \begin{eqnarray*}
  |w'|&=&\left|\sum_{k=0}^{p -1}{\bar{z'}^k}A_k\left( z' \right)\right|\geq |A_0(z')| - \left|\sum_{k=1}^{p -1}{a_{k,1}}\left| z \right|^2\bar{z'}^{k-1}+\sum_{k=1}^{p -1}{\bar{z'}^k\sum_{n=2}^{\infty}{a_{k,n}z'^n}} \right|\\
  &\geq&\Lambda^2 \rho _4+(\Lambda^3-\Lambda)\ln\left(1-\frac{\rho _4}{\Lambda}\right)-\sum_{k=1}^{p -1}{\rho _{4}^{k+1}}-\sum_{k=0}^{p -1}{\left( M_k-\frac{1}{M_k} \right)}\sum_{n=2}^{\infty}{\rho _{4}^{n+k}}\\
  &=&\Lambda^2 \rho _4+(\Lambda^3-\Lambda)\ln\left(1-\frac{\rho _4}{\Lambda}\right)-\rho _{4}^{2}\left( \frac{1-\rho _{4}^{p -1}}{1-\rho _4} \right) -\sum_{k=1}^{p -1}{\left( M_k-\frac{1}{M_k} \right)}\frac{\rho _{4}^{k+2}}{1-\rho _4}=\sigma_4,
  \end{eqnarray*}
  which implies that $F(U_{\rho_3})\supseteq U_{\sigma_3}$. \hfill $\Box$

 \subsection{Proof of Theorem \ref{HLP-th5}}
 Let $F(z)=\sum\limits_{k=0}^{p-1}{\bar{z}^k}A_k\left( z \right)$. Then it follows from the hypothesis of Theorem \ref{HLP-th5} and the definition of log-analytic function that for each $k\in \left\{ 0,1,\cdots ,p-1 \right\}$, $A_k(z)=\log a_k(z)$ is an analytic function in $U$. Thus $F(z)=\log f(z)$ is a $p$-analytic function in $U$.

 We first prove the univalence of $f(z)$ in $U_{\rho_1}$. Indeed, for any $z_1,z_2\in U_r\, (0<r<\rho_1$, $z_1\neq z_2)$. Then it follows from our proof of Theorem \ref{HLP-th1} and the hypothesis of Theorem \ref{HLP-th5} that
 \begin{eqnarray*}
|\log f\left( z_1 \right) -\log f\left( z_2 \right) |&=&|F(z_1)-F(z_2)|
=\left| \int_{\left[ z_1,z_2 \right]}{F_z\left( z \right)}dz+F_{\bar{z}}\left( z \right) d\bar{z} \right|\nonumber\\
 &\geq&|z_1-z_2|\left(\Lambda_0 \frac{1-\Lambda_0 r}{\Lambda_0-r}-\sum_{k=1}^{p -1}{\left( k+1 \right)}\varLambda _kr^k\right).\nonumber
 \end{eqnarray*}

 Thus $\log f\left( z_1 \right) \ne \log f\left( z_2 \right) $, which implies the univalency of $f$ in the disk $U_{\rho_1}$.

 Next, we consider any $z'$ with $|z'|=\rho_1$. By our proof of Theorem \ref{HLP-th1}, we have that
 \begin{eqnarray*}
 |\log f\left( z' \right)|&=&|F(z')|=\left|\sum_{k=0}^{p -1}{\bar{z'}^k}A_k\left( z' \right)\right|\geq\left| A_0\left( z' \right) \right|-\sum_{k=1}^{p -1}{\rho _{1}^{k}\left| A_k\left( z' \right) \right|}\\
 &\geq&\Lambda_0^2 \rho_1-\sum_{k=1}^{p -1}{\varLambda _k}\rho_1^{k+1}+\left(\Lambda_0^3-\Lambda_0\right)\ln\bigg(1-\frac{\rho_1}{\Lambda_0}\bigg)=\sigma_1.
 \end{eqnarray*}
 where $\sigma_1$ is defined by \eqref{liu22}.

 By Lemma \ref{HLP-lem3}, we obtain that the range $f\left( U_{\rho_1} \right)$ contains a univalent disk $U\left( w_1,r_1 \right) $, where $w_1$ and $r_1$ are defined by \eqref{liu210}.

 Now, we prove the sharpness of $\rho_1$ and $\sigma_1$. To this end, we consider the log-p-analytic function $f_1(z)=e^{F_1(z)}$, where $F_1(z)$ is given by \eqref{liu23}. It is easy to verify that $f_1(z)$ satisfies the hypothesis of Theorem \ref{HLP-th5},
and thus, we have that $f_1(z)$ is univalent in $U_{\rho_1}$, and the range $f_1\left( U_{\rho_1} \right) $ contains a schlicht disk $U\left( w_1,r_1 \right)$.

 To prove that the univalent radius $\rho_1$ is sharp, we need to prove that $f_1(z)$ is not univalent in $U_r$ for each $r\in (\rho_1, 1]$. In fact,
if we fix $r\in (\rho_1, 1]$, by our proof of Theorem  \ref{HLP-th1}, we know that $F_1(z)$ is not univalent in $U_r$, thus there exist  two distinct points $z_1,z_2 \in U_r$ such that $F_1(z_1)=F_1(z_2)$, which implies that $$f_1(z_1)=e^{F_1(z_1)}=e^{F_1(z_2)}=f_1(z_2)$$
in that, $f_1(z)$ is not univalent in $U_r$ for each $r\in (\rho_1, 1]$. Hence, the univalent radius $\rho_1$ is sharp.

Finally, we prove that radius $r_1=\sinh\sigma_1$ is sharp. In fact, by the proof of Theorem  \ref{HLP-th1}, we have that $\sigma_1>0$ and  $0<\rho_1\leq 1/\varLambda _0$. We now prove that $\sigma_1<1$. Let
$$
h(r)=\Lambda_{0}^{2}r+\left(\Lambda_{0}^{3}-\Lambda_{0}\right)\ln\left(1-\frac{r}{\Lambda_{0}}\right),\quad 0<r\leq\frac{1}{\varLambda _0},
$$
then
$$
h'(r)=\Lambda_0^2+\frac{1-\Lambda_0^2}{1-\frac{r}{\Lambda_0}}=\Lambda_0\frac{\frac{1}{\Lambda_0}-r}{1-\frac{r}{\Lambda_0}}\ge0,\quad 0<r\le\frac{1}{\varLambda _0},
$$
which implies that $h(r)$ is increasing in $(0,1/\Lambda_0]$, therefore
\begin{eqnarray*}
 \sigma_1&=&\Lambda_0^2 \rho_1-\sum_{k=1}^{p -1}{\varLambda _k}\rho_1^{k+1}+\left(\Lambda_0^3-\Lambda_0\right)\ln\bigg(1-\frac{\rho_1}{\Lambda_0}\bigg)\\
 &\leq& h(\rho_{1})\leq h(\frac{1}\Lambda_0)=\Lambda_{0}+\left(\Lambda_{0}^{3}-\Lambda_{0}\right)\ln\left(1-\frac{1}{\Lambda_{0}^{2}}\right) \\
 &<&\Lambda_0+\left(\Lambda_0^3-\Lambda_0\right)\cdot\left(-\frac{1}{\Lambda_0^2}\right)=\frac{1}{\Lambda_0}<1.
\end{eqnarray*}
Hence, $0<\sigma_1<1$.

Given that the univalent radius \( \rho_1 \) is exact, the precision of the radius \( r_1 = \sinh(\sigma_1) \) is inferred from Lemma \ref{HLP-lem3}, considering the condition \( 0 < \sigma_1 < 1 \). With this, the demonstration of the theorem is concluded. \hfill $\Box$

\subsection{Proof of Theorem 7}
By the hypothesis of Theorem 7, we have that $A_k(z)=\log a_k(z)$ is analytic for each $k\in \left\{ 0,1,\cdots ,p-1 \right\}$, which implies that
 $$
 F(z)=\sum\limits_{k=0}^{p-1}{\bar{z}^k}A_k\left( z \right)
 $$
 is a $p$-analytic function in $U$.

 We first prove the univalence of $f(z)$ in $U_{\rho_3}$. For $k\in \left\{ 0,1,\cdots ,p-1 \right\}$, $z\in U$, we have
 $$
 |A_k(z)|=|\log a_k(z)|=|\log|a_k(z)|+i\arg a_k(z)|\le |\log|a_k(z)||+\pi,
 $$
 This implies that $|A_k(z)|\le\log M_{k}^{*}+\pi:=M_k$. Therefore, for any $z_1,z_2\in U_r\, (0<r<\rho_3$, $z_1\neq z_2)$. It follows from our proof of Theorem \ref{HLP-th3} and the hypothesis of Theorem 7 that
 \begin{align*}
&|\log f\left( z_1 \right) -\log f\left( z_2 \right) |=|F(z_2)-F(z_1)|=\left| \int_{\left[ z_1,z_2 \right]}{F_z\left( z \right)}dz+F_{\bar{z}}\left( z \right) d\bar{z} \right|\nonumber\\
 &\geq|z_{2}-z_{1}|\left(1-{\sum_{k=0}^{p-1}}\left( M_k-\frac{1}{M_k} \right)\frac{r^{k+1}(2-r+k(1-r))}{(1-r)^{2}}-{\sum_{k=1}^{p-1}}(k+1)r^{k}\right)>0.\nonumber
 \end{align*}

 Thus $\log f\left( z_1 \right) \ne \log f\left( z_2 \right) $,which implies the univalency of $f$ in the disk $U_{\rho_3}$.

 Next, we consider any $z'$ with $|z'|=\rho_3$. By our proof of Theorem \ref{HLP-th3}, we have
 \begin{eqnarray*}
 |\log f\left( z' \right)|&=&|F(z')|=\left|\sum_{k=0}^{p -1}{\bar{z'}^k}A_k\left( z' \right)\right|\\
 &\geq&\rho _3-\sum_{k=1}^{p -1}{\rho _{3}^{k+1}}-\sum_{k=0}^{p -1}{\left( M_k-\frac{1}{M_k} \right)}\sum_{n=2}^{\infty}{\rho _{3}^{n+k}}\\
&=&\rho_3-\rho _{3}^{2}\left( \frac{1-\rho _{3}^{p -1}}{1-\rho _3} \right) -\sum_{k=0}^{p -1}{\left( M_k-\frac{1}{M_k} \right)}\frac{\rho _{3}^{k+2}}{1-\rho _3}=\sigma_3.
 \end{eqnarray*}
 where $\sigma_3$ is defined by \eqref{liu29}.

 By Lemma \ref{HLP-lem3}, we obtain that $f\left( U_{\rho_3} \right) $ contains a schlicht disk $U\left( w_3,r_3 \right) $, where $w_3$ and $r_3$ are defined by \eqref{liu212}.
The proof of this theorem is complete.\hfill $\Box$

\section{Conclusion}\label{HLP-sec5}

In conclusion, this paper makes a significant contribution to the field of complex analysis by introducing the innovative concept of log-$p$-analytic functions. The study meticulously develops four unique Landau-type theorems for poly-analytic functions, with two of these theorems being particularly sharp and another providing a notable refinement to the work of Abdulhadi and Hajj. The paper culminates in the presentation of four specialized Landau-type theorems for bounded log-$p$-analytic functions, yielding two exact results that advance the understanding of these mathematical constructs. This research not only enriches the theoretical framework but also opens avenues for further exploration and application within the domain of complex analysis.

\end{document}